\documentclass[11pt]{amsart}

\usepackage{amssymb, graphicx}

\newcommand{\ncmnd}{\newcommand}
\ncmnd{\nthm}{\newtheorem}
\numberwithin{equation}{section}

\theoremstyle{plain}
\nthm{thm}{Theorem}[section]
\nthm{cor}[thm]{Corollary}
\nthm{prop}[thm]{Proposition}
\nthm{lem}[thm]{Lemma}

\theoremstyle{definition}
\nthm{defn}[thm]{Definition}
\nthm{examp}[thm]{Example}

\theoremstyle{remark}
\nthm{rem}[thm]{Remark}

\newcommand{\R}{\mathbf{R}}
\newcommand{\E}{\mathbf{E}}

\begin{document}

\title{The intrinsic geometry of a Jordan domain}
\author{Richard L. Bishop}
\address{1409 W. Green St., Urbana, Illinois 61801}
\email{bishop@math.uiuc.edu}

\begin{abstract}
For a Jordan domain in the plane  the length metric space of points connected
to an interior point by a curve of finite length is a {\sc cat(0)} space and
Gromov hyperbolic. With respect to the cone topology, that space plus
its boundary at infinity is topologically the same as the original Jordan domain.
\end{abstract}

\maketitle

\section{Introduction}
\label{sec_intro}
A Jordan domain $J$ is the the homeomorphic image of a closed disk in
$\E^2$. The image of the boundary circle is a Jordan curve, which by
the Jordan Curve Theorem separates the plane into two open domains,
one bounded, the other not, such that the curve is the boundary of
each. A {\em local geodesic} is a curve $\gamma$ such that for every non-end
point $p$ of $\gamma$ there is an open subarc $\rho$ of $\gamma$
containing $p$ such that the arc of $\rho$ between any two points of
$\rho$ is a shortest curve in $J$ connecting those two points.
Bourgin and Renz \cite{BR89} have analyzed the local geodesics
in such domains, concluding the following:
\begin{itemize}
\item A Jordan domain is uniquely geodesic.
\item At every non-end point of a geodesic, the geodesic is supported by a
closed half-disk with center at the point and interior contained in
the interior of the Jordan domain.
\item Geodesics are characterized in terms of separation: a point
$p$ is not on the geodesic from $q$ to $r$ if and only if there is a
straight line segment $\sigma$ with ends on the boundary $\partial J$ of $J$,
but otherwise in the interior of $J$, such that $p$ is in one connected
component of $J \setminus \sigma$ and $q,r$ are in the other.
\end{itemize}
Although geodesics may be infinitely long when an end is on the boundary, the
unique determination of a geodesic by its ends is still true.

We assume throughout that geodesics are parametrized by an arclength
parameter, and for a geodesic ray the parameter is $0$ at the origin of the
ray and positive elsewhere.

There has been a little previous
related work concerning not quite the same domains: in the thesis of
F.-E. Wolter \cite{W85} it is proved that a simply connected plane domain
with locally rectifiable boundary has the unique geodesic property; in a recent
text \cite[p. 310]{BBI01} there is an outline of a proof that locally
simply-connected plane domains are locally {\sc cat(0)}. In the latter,
although the outline is entirely plausible and at a level of rigor appropriate
to a text at that level, essential technical details (the unique geodesic
property and a consideration of triangles with ``tails'') are not mentioned.
It is probable that the {\sc cat(0)} property, but certainly not Gromov
hyperboliciy, can be extended to unbounded domains by exhausting them
with their intersection with increasing disks; also extensions giving natural
curvature bounds for similar domains in Riemannian surfaces are likely to
be true; but it seems as though these extensions would not present any
interesting challenges.

\section{The {\sc cat(0)} space}

The set $X$ of {\em finite-distance} points of $J$ consists of interior
points and those points of $\partial J$ which can be connected to any
(and hence every) interior point by a curve of finite length. Note that
every pair of points in the interior can be connected by polygonal curve
in the interior, and so have finite distance apart. If a straight
line segment is interior to $J$ except for one end, then that end is also
in $X$; in particular, only the ends of a geodesic can fail to be in $X$,
since non-ends are either already interior to $J$ or a radius of a supporting
half-disk provides a segment to an interior point. The remaining subset
of $\partial J$, points not in $X$, is called the {\em (metric) boundary}
of $X$, denoted $\partial X$.

A {\em geodesic triangle} in a metric space consists of three shortest curves
(its {\em sides}) connecting a triple of points (its {\em vertices}).
A {\em comparison triangle} in $\E^2$ is a triangle with the same
side-lengths, and by taking equal distances from the vertices we get
comparison points on the sides. A geodesic triangle is {\em thin} if the
distance between every pair of points in the sides is $\le$ the distance
between the comparison points. A {\sc cat(0)} {\em space} is a complete
metric space such that every pair of points can be joined by a (finite-length)
geodesic and every geodesic triangle is thin. It is easy to show that a
{\sc cat(0)} space is uniquely geodesic and contractible \cite{BH99}.
Conversely, if a space is locally {\sc cat(0)} and simply connected,
then it is {\sc cat(0)} \cite{AB90}.
\begin{thm} \label{thm:main} The set $X$ of finite-distance points of a Jordan
domain $J$ is a {\sc cat(0)} space.
\end{thm}
\begin{proof} Consider three points $p,q,r \in X$. We suppose that
the geodesic triangle with vertices $p,q,r$ is nondegenerate, since
degenerate triangles are trivially thin. Then the two sides starting
from $p$ must consist of a common part (which may be just $p$) ending
at a bifurcation point $\bar p$; similarly, there are bifurcation
points $\bar q, \bar r$ on the sides starting from $q,r$. Then
$\triangle\bar p \bar q \bar r$ is a Jordan curve in $X$ and
$\triangle pqr$ consists of $\triangle\bar p \bar q \bar r$ with three
``tails''. Clearly if we show that $\triangle\bar p \bar q \bar r$ is
thin, then $\triangle pqr$ is also thin.

The interior of $\triangle\bar p \bar q \bar r$ must lie in the
interior of $X$, and at each point of a side which also lies in the
interior of $X$, that side is locally a straight line segment. Hence
at those points the side is supported by half-disks in the interior of
$X$, both interior and exterior to $\triangle\bar p \bar q \bar
r$. But at a point of a side which is in $\partial X$, a supporting
half-disk can only lie entirely in the domain bounded by
$\triangle\bar p \bar q \bar r$. Thus, a side of $\triangle\bar p \bar
q \bar r$ is supported at every point except the ends by an open
half-disk interior to $\triangle\bar p \bar q \bar r$. Therefore the
sides are locally convex and have well-defined one-sided directions at
every point and also a total curvature. At the ends the signed turning
angle (choose an orientation of $\E^2$!) is at most $\pi$, while the
signed total curvature of any simple closed, locally convex curve must
be $2\pi$. The signed total curvatures of the sides must be $\le 0$,
so together they add to at least $2\pi - 3\pi = -\pi$. In
particular, each side must actually be convex, i.e., is on the
boundary of its convex hull, and the remaining part of that boundary
is a straight line segment. The triangle of these three segments has
$\triangle\bar p \bar q \bar r$ in its interior, so that
$\triangle\bar p \bar q \bar r$ must have angle sum $\le 2\pi$.  One
of Alexandrov's criteria for a locally {\sc cat(0)} space is that locally angle sums
of triangles are $\le 2\pi$. Since $X$ is simply connected, we
conclude that it is also {\sc cat(0)}.

It is also easy to see that $X$ is complete. Indeed, distances in $X$
are $\ge$ the corresponding Euclidean distances, so that a Cauchy
sequence in $X$ is also Cauchy in the plane. The limit in the plane is
in $J$ and the distance of the limit from any point of the sequence is
finite.
\end{proof}

Let us call a triangle in $X$ such as $\triangle\bar p \bar q \bar r$
a {\em Jordan triangle}. As a corollary to the above proof we have a
uniform bound on the perimeter of Jordan triangles:
\begin{cor} \label{cor:trianglebound} If $D$ is the diameter of
$J$ in $\E^2$, then the perimeter of every Jordan triangle is
bounded by $4D$.
\end{cor}
\begin{proof} In $\E^2$ the length of a closed convex curve is
monotone increasing with respect to the inclusion ordering of convex
hulls. We can choose Euclidean segments from the vertices of
$\triangle\bar p \bar q \bar r$ so that they don't intersect the sides and
either are concurrent or meet in pairs to give a triangle $\triangle abc$ 
in the interior. The latter case is generic, pictured in Figure \ref{fig:Jtriangle}.
The sides of the Jordan triangle are shorter than the smaller triangle in
which they are contained; that is, $\ell(\tau) < d(\bar p, a) + d(a, \bar r)$, etc.
In the concurrent case, when point $a = b = c = \bar m$, then
the segments $\bar m \bar p, \,\bar m \bar q, \,\bar m \bar r$
are in the interior of $\triangle\bar p \bar q \bar r$
except for the ends. The sum of the lengths of these segments from
$\bar m$ is at most $2D$ (the sum of distances in $\E^2$ from
several points is a convex function, so the maximum on a triangle must
occur at a vertex), and hence twice that sum is greater than twice the
sum of the two longest side-lengths of $\triangle\bar p \bar q \bar r$.

\begin{figure}[b]
\begin{center}
\includegraphics[width= 3in]{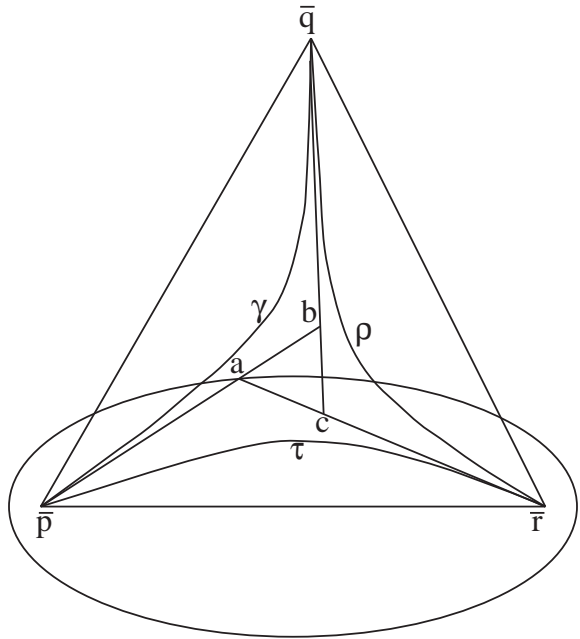}
\end{center}
\caption{Jordan triangle}
\label{fig:Jtriangle}
\end{figure}

The generic case can be reduced to the concurrent case by showing that
there is an interior point $\bar m$ of $\triangle abc$ such that the segments
$\bar m \bar p, \,\bar m \bar q, \,\bar m \bar r$ again have pairwise
sums of lengths which bound the lengths of $\tau, \gamma, \rho$.
Consider the ellipse with foci $\bar p, \bar r$ passing through $a$. For
any point on that ellipse the sum of distances from $\bar p, \bar r$ is
the same, hence greater than $\ell(\tau)$. For a point outside the ellipse
the sum is even greater, so that is true for any point (except $a$) on the
tangent line to the ellipse at $a$. But the rule for reflection in an ellipse
tells us that that tangent line bisects $\angle bac$. Similarly, the angle
bisectors of $\angle abc, \angle bca$ give points from which the sums
of distances from $\bar p,\bar q$ and $\bar q, \bar r$ majorize the
lengths of $\gamma, \rho$. Thus we can take $\bar m$ to be the incenter
of $\triangle abc$.
\end{proof}

\section{Gromov hyperbolicity}
A geodesic metric space $X$ is {\em$\delta$-hyperbolic}, where
$\delta \ge 0$, if for every triangle the distance from any point on one
side of the triangle to the union of the other two sides is $\le \delta$.
If $X$ is $\delta$-hyperbolic for some $\delta$, then it is
{\em Gromov hyperbolic}.

When $X$ is a Jordan domain, it is clear that it is $\delta$-hyperbolic if the defining property is true for every Jordan triangle in $X$. Half the perimeter of the triangle obviously satisfies the defining property, so by Corollary \ref{cor:trianglebound} $X$ is $2D$-hyperbolic. However, a better value for $\delta$ is obtained by observing that for a Euclidean triangle the extreme case is equilateral and if the side-lengths are $D$, then the maximum distance from a point to the other sides is half the altitude. Thus we obtain the sharp value for $\delta$, realized by the case when $X$ is bounded by an equilateral triangle of side $D$.
\begin{cor}If $X$ is a Jordan domain with Euclidean diameter $D$, then it is $\sqrt{3}D/4$-hyperbolic.
\end{cor}

\section{The cone topology}

For a {\sc cat(0)} space $X$ the boundary at infinity $\partial X$
consists of asymptote classes of geodesic rays; two rays are {\em
asymptotic} if the distance between pairs at equal distance from their
origins is bounded. The {\em cone topology} on $\bar X = X\cup\partial
X$ is defined by specifying a neighborhood basis: for points of $X$ we
take the usual metric neighborhoods. The neighborhoods of $\partial X$
are defined in terms of a fixed center point $p \in X$, and then it is
proved that the resulting topology is independent of $p$. For an
infinite ray $\gamma$ starting at $p$ and positive numbers $C,
\epsilon$ we define a neighborhood $\mathcal N(\gamma, C,\epsilon)$ of
the asymptote class of $\gamma$ to be the points on geodesics starting
at $p$ which pass through the ball $B(\gamma(C), \epsilon)$ and have
distance $>C$ from $p$, along with the asymptote classes of the
extensions of those geodesics to rays whenever possible. (The
definition was originally due to Eberlein and O'Neill for Riemannian
manifolds \cite{EO73}.) For each asymptote class there is exactly one
ray for a given origin \cite{BH99}.

\begin{lem} \label{lem:rays}For a Jordan domain $J$ and the space of
finite points $X$, two asymptotic rays in $X$ must eventually
coincide.
\end{lem}
\begin{proof} We first prove this for two asymptotic rays
$\gamma, \sigma$ for which the origins are interior points $p,q$ of
$J$ such that the geodesic connecting them is a straight line segment
in the interior. Then any triangle $\triangle pqr$ must end in a tail
at $r$ as soon as $d(p,r) + d(q,r) > 4D$, since bifurcation of
geodesics can only occur at boundary points, not at $p,q$ (unless
$\triangle pqr$ is degenerate and the following conclusion is
trivial). Taking $r = \gamma(s)$, as soon as $s + d(q,r) > 4D$ the
continuation of the segment from $q$ to $\gamma(s)$ by
$\gamma|_{[s,\infty)}$ is a geodesic ray asymptotic to $\gamma$. Since
the origin is $q$, that ray must coincide with $\sigma$. So for
close-by interior points ``eventually'' means as soon as the sum of
distances from their origins exceeds $4D$.

Now for arbitrary origins in $X$ we can connect them by a curve $\tau$
through the interior, and for any two interior points of $\tau$,
we conclude that beyond the distance $2D$ from
$\tau$, the rays in the given asymptote class coincide. The
limit rays at the ends then must also coincide beyond that distance.
\end{proof}

\begin{thm}\label{thm:topology} For a Jordan domain $J$ the asymptote
classes of the space of finite points $X$ can be indentified with the
points of $\partial J$ at infinite distance. The topology on $J$
coincides under this identification with the cone topology on $X\cup\partial X$.
\end{thm}
\begin{proof} A ray in $X$ has a unique limit point in $\partial J$,
which by Lemma \ref{lem:rays} is the same for all rays in the same
asymptote class.

Suppose that $q \in \partial J$ is the end of a ray $\gamma$ from
$p \in X$, so $q$ is identified with the asymptote class of $\gamma$.

Let $B(\gamma(C), \epsilon)$ specify a cone neighborhood
$\mathcal N(\gamma, C,\epsilon)$ of $q$. Then there is a straight line
segment $rs$ contained in $B(\gamma(C), \epsilon)$, and
containing $\gamma(C)$, transverse to $\gamma$. If
$\gamma(C) \in \partial J$, then, say, take $r = \gamma(C)$ and
$rs$ to be a radius of a supporting half-disk at $\gamma(C)$. Otherwise
$rs$ can be in the interior of $J$. Let $m, n$ be the first points in
geodesic extensions of the geodesics $pr, ps$ where those  extensions
intersect $\partial J$ (so that if $r \in \partial J$ let $m = r$).
Since $m, n$ must be at finite distance from $p$, they are different
from $q$. See Figure \ref{fig:cone}.

\begin{figure}[b]
\begin{center}
\includegraphics[width= 4in]{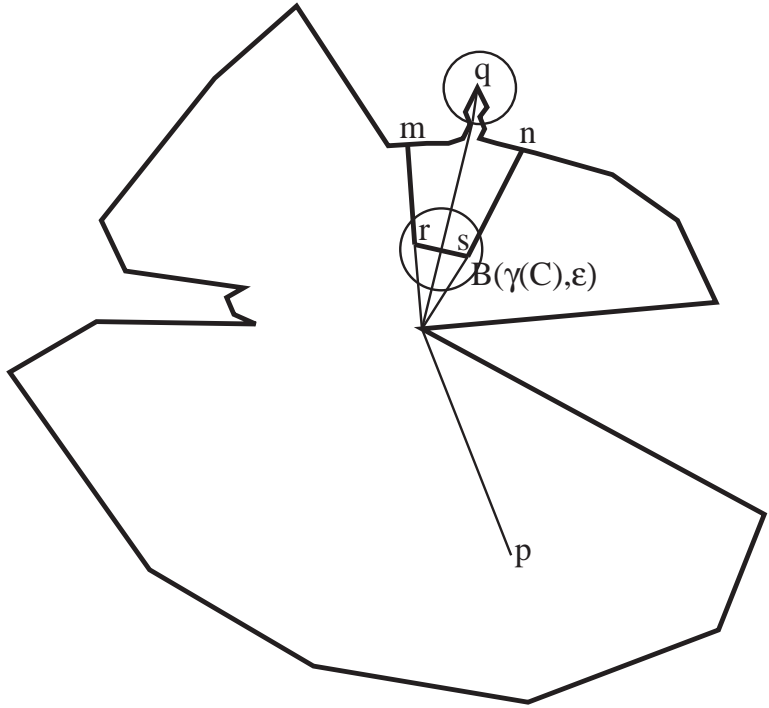}
\end{center}
\caption{Disk inside cone neighborhood}
\label{fig:cone}
\end{figure}

Then there is a Jordan curve obtained
by chaining together geodesics $mr$, $rs$, $sn$ with the arc of
$\partial J$ from $n$ to $m$ which does not contain $q$. A sufficiently
small disk centered at $q$ will not intersect that new Jordan curve, and
any geodesic from $p$ to a point of that disk must cross $rs$, so the
point in the disk must be in $\mathcal N(\gamma, C,\epsilon)$.
Thus there are neighborhoods of the induced topology on $J$ in a
cone topology neighborhood.

Conversely, suppose that $U$ is the intersection of $J$ with a disk $D_1$ centered
at $q$, so that $U$ is a neighborhood of the induced topology. We may
assume that $U$ does not contain $p$. Then there
is a largest arc of $\partial J$ containing $q$ and contained in $U$. Let
$m, n$ be the ends of that arc. Since $q$ is at infinite distance, it cannot
be on the geodesic from $m$ to $n$. Hence there must be a straight line
segment with ends $r, s \in \partial J$ separating $mn$ from $q$. See
Figure \ref{fig:disk}. 
\begin{figure}[b]
\begin{center}
\includegraphics[width= 4in]{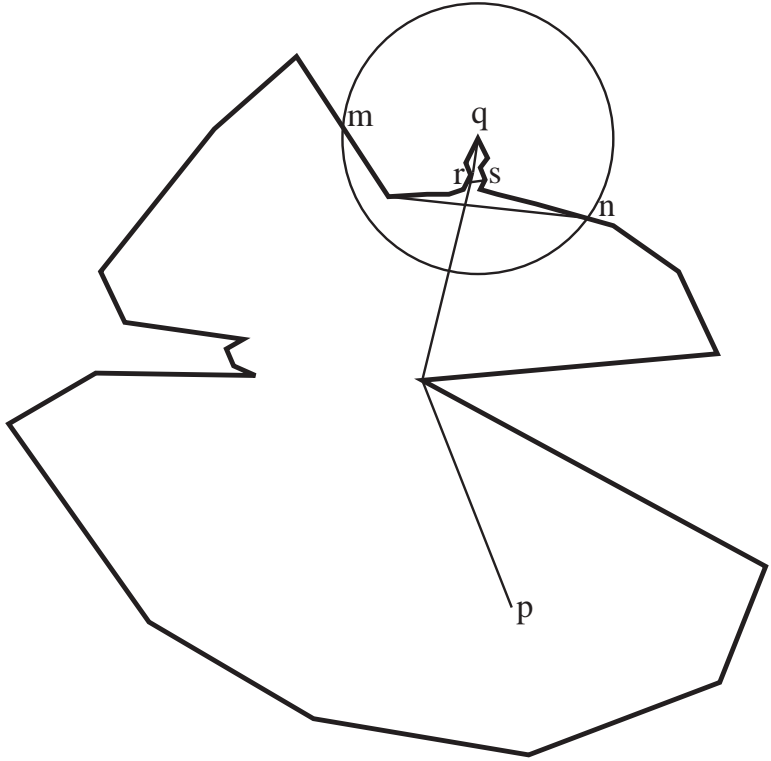}
\end{center}
\caption{Cone neighborhood inside disk}
\label{fig:disk}
\end{figure}
\noindent The arcs of $\partial J$ in $U$ from $m$ to $q$
and from $n$ to $q$ each must contain only one of $r, s$, since both
of those arcs connect $q$ to $\{m,n\}$ and cross $rs$. The geodesic
$mn$ together with the arc of $\partial J$ with ends $m, n$ not
containing $q$ is a Jordan curve having $p$ inside and $q, rs$ outside;
then the Jordan curve consisting of $rs$ and the arc of $\partial J$
from $r$ to $s$ through $q$ must bound a domain $D_2$ inside $U$. Any
geodesic from $p$ to a point of $D_2$ must cross $rs$. Let
$\gamma(C)$ be the point where $pq$ crosses $rs$. Then for sufficiently
small $\epsilon$ the ball $B(\gamma(C),\epsilon)$ will be separated by
$rs$ and contained in $U$; we can take $\epsilon <$ the distance in $\E^2$
from $\gamma(C)$ to $\partial D_1$. Then any point of
$\mathcal N(\gamma, C,\epsilon)$ will either be in $B(\gamma(C),\epsilon)$
or $D_2$, and hence in $U$.
\end{proof}
\begin{rem} The continuity of geodesics in $J$ as a function of their
endpoints has also been addressed by Fabel (\cite{F99}). He proves that
when the endpoints of a sequence converge in $J$, then the geodesics
converge uniformly with respect to the Euclidean metric. There is also
a convergence interpretation of Theorem \ref{thm:topology}. Namely,
it says that the convergence of endpoints is equivalent to convergence
of the geodesics uniformly on bounded sets in the length metric. A related
fact for {\sc cat(0)} spaces, which provides a link between the two kinds
of convergence in Jordan domains, is that the distance
$d(\gamma(t), \sigma(at))$
between two geodesics is a convex function of $t$ for any $a>0$. On a
bounded interval of $\R$ the convergence of the end values of a convex
function implies the uniform convergence on the interval. Thus Euclidean
convergence of the infinite endpoints controls the unbounded ends of geodesics
uniformly with respect to the Euclidean metric, and the convexity of the
length metric (which majorizes the Euclidean metric) controls the bounded
part. This provides the kernel of the proof that the two kinds of convergence
are equivalent. However, Fabel goes on to prove a much stronger theorem,
namely, the Euclidean-uniform continuity of geodesics not only as functions
of the ends, but also as functions of the domain $J$ with respect to the
Hausdorff distance.
\end{rem}

\section*{Acknowledgments} This paper was an outgrowth of another
project by S. Alexander, R. Ghrist, and the author concerning the use of
{\sc cat(0)} spaces in pursue-evade theory. In particular, the author
thanks Ghrist for passing on the reference \cite{F99} to him.

\end{document}